\documentclass[12pt, a4paper]{article}
\usepackage{amsfonts}
\usepackage{mathrsfs}

\usepackage{latexsym}
\usepackage{graphicx}
\usepackage{amssymb}
\usepackage{bm}
\usepackage{amsmath}

\newtheorem{theorem}{Theorem}[section]
\newtheorem{lemma}[theorem]{Lemma}

\newtheorem{definition}[theorem]{Definition}

\title{{\bf Planar Tur\'{a}n  Number of intersecting triangles}
\thanks{Supported by the National Natural Science Foundation of China (Nos. 41571398 and 11871222) and SRSF of Chuzhou University (No. 2018qd02).}}

\author{{\bf Longfei Fang}, {\bf Mingqing Zhai}\thanks{Corresponding author: mqzhai@chzu.edu.cn
(M.Zhai)}, {\bf Bing Wang}
\\
{\footnotesize School of Mathematics and Finance, Chuzhou
University, Anhui, Chuzhou, 239012, China}}
\date{}

\begin{document}
\openup 1.0\jot \maketitle

\begin{abstract}
The planar Tur\'{a}n number of a given graph $H$, denoted by $ex_{\mathcal{P}}(n,H)$, is the maximum number of edges over all planar graphs on $n$ vertices that do not contain a copy of $H$ as a subgraph. Let $H_k$ be a friendship graph, which is obtained from $k$ triangles by sharing a common vertex. In this paper, we obtain sharp bounds of $ex_{\mathcal{P}}(n,H_k)$ and $ex_{\mathcal{P}}(n,K_1+P_{k+1})$ for $k\ge2$, which improve the results of Lan, Shi and Song in Electron. J. Combin. 26 (2) (2019), \#P2.11.


\end{abstract}

\section{Introduction}

~~~~Given a graph $H$, a graph is said to be \emph{$H$-free} if it contains no $H$ as a subgraph.
As one of the best results in extremal graph theory, Tur\'{a}n Theorem gives
the maximum number of edges in a $K_r$-free graph on $n$ vertices. Replacing $K_r$ with an arbitrary graph $H$, the famous Erd\H{o}s-Stone Theorem shows that the maximum number of edges in an $H$-free graph on $n$ vertices possible is $(1+o(1)){{n}\choose{2}}\Big(\frac{\chi(H)-2}{\chi(H)-1}\Big)$, where $\chi(H)$ denotes the chromatic number of $H$. The \emph{Tur\'{a}n number} of a graph $H$, denoted by $ex(n, H)$, is the maximum number of edges in a graph $G$ on $n$ vertices which does not contain $H$ as a subgraph. Over the last decades, a large quantity of work has been carried out in
Tur\'{a}n-type problems (see the survey paper \cite{KB}).

In 2015, Dowden \cite{DZ} initiated the study of Tur\'{a}n-type problems when the host graph is planar, i.e., how many edges can an $H$-free planar graph on $n$ vertices have? The \emph{planar Tur\'{a}n number} of a graph $H$, denoted by $ex_{\mathcal{P}}(n,H)$, is the maximum number of edges in an $H$-free planar graph on $n$ vertices. It is easy to check that a triangle-free planar graph has at most $2n-4$ edges and $K_{2,n-2}$ is exactly the extremal graph. Besides, there is a plane triangulation containing no $K_4$ (see, for example, the graph $2K_1+C_{n-2}$).
Thus, $ex_{\mathcal{P}}(n,K_3)=2n-4$ and $ex_{\mathcal{P}}(n,K_r)=3n-6$ for all $r\ge 4$. Since the planar Tur\'{a}n problem on complete graphs is trivial, the next natural type of graphs considered are cycles and others variations. Dowden \cite{DZ} showed that $ex_{\mathcal{P}}(n,C_4)\le \frac{15(n-2)}{7}$ for  $n\ge 4$ and $ex_{\mathcal{P}}(n,C_5)\le \frac{12n-33}{5}$ for $n\ge 11$. Ghosh et. al \cite{DG} proved that $ex_{\mathcal{P}}(n,C_6)\le \frac{5n}{2}-7$ for $n\ge 18$. Lan et. al \cite{BT} obtained a sharp upper bound of $ex_{\mathcal{P}}(n,\mathcal{C}_k+e)$ for $k\in \{4,5\}$ and gave an upper bound of $ex_{\mathcal{P}}(n,\mathcal{C}_6+e)$, where $\mathcal{C}_k+e$ is the family of graphs obtained from a cycle $C_k$ by linking two nonadjacent vertices (via $e$) in the cycle. Later, Ghosh et. al \cite{LN} obtained a sharp bound of $ex_{\mathcal{P}}(n,\mathcal{C}_6+e)$, which improved the result of Lan et. al.

In this paper, we focus on the planar Tur\'{a}n number of two classes of  graphs.
Let $H_k$ be a \emph{friendship graph}, which is obtained from $k$ triangles by sharing a common vertex. Let $F_k\cong K_1+P_{k+1}$ be a \emph{$k$-fan}, which is obtained by joining a new vertex to all the vertices of $P_{k+1}$. In 1995, Erd\H{o}s et. al \cite{DX} determined $ex(n,H_k)$ for every fixed $k$ and large enough $n$. Chen et. al \cite{DW} studied $ex(n,K_1+kK_r)$ for sufficiently large $n$, where $k\ge 2$ and $r\ge 2$. Recently, Yuan \cite{YLT} determined $ex(n,F_k)$ for $k\ge 3$ and large enough $n$. Now we consider the case that the host graphs are planar. Notice that the average degree of a plane triangulation is less than $6$ and for each $n$, there are plane triangulations of order $n$ such that the maximum degree are at most $6$. Therefore, $ex_{\mathcal{P}}(n,H_k)=3n-6$ for $k\ge 4$ and $ex_{\mathcal{P}}(n,F_k)=3n-6$ for $k\ge 6$. Lan, Shi and Song \cite{D} obtained a sharp upper bound of $ex_{\mathcal{P}}(n,H_2)$ for $n\ge 5$, and then they gave the following result for $ex_{\mathcal{P}}(n,H_3)$.

\begin{theorem}\label{503}(\cite{D})
	Let $n\ge 7$ be an integer. Then
	$$\lfloor\frac{5n}{2}\rfloor \le ex_{\mathcal{P}}(n,H_3)<\frac{17n}{6}-4$$
	for all $n\ge 15$ and
	
	$$  ex_{\mathcal{P}}(n,H_3)=\left\{
	\begin{array}{rcl}
	3n-6     &    & {if~n\in\{7,8,9,10,12\}}\\
	3n-7     &    & {if~n=11}\\
	3n-8     &    & {if~n\in\{13,14\}}
	\end{array} \right. $$
\end{theorem}

However, the upper bound in Theorem \ref{503} is not sharp. In Section 2, we show a sharp upper bound of $ex_{\mathcal{P}}(n,H_3)$ and give infinitely many extremal graphs which attain the extremal values.

  Now, let us focus our attention on $k$-fans.
Recall that $ex_{\mathcal{P}}(n,F_k)=3n-6$ for $k\ge 6$.
 Lan, Shi and Song also considered the values of $ex_{\mathcal{P}}(n,F_k)$ for the remaining cases.

\begin{theorem}\label{103}
Let $n\ge k+2$. Then\\
(\romannumeral1) (\cite{BT}) for $k=2$, $ex_{\mathcal{P}}(n,F_k)\le \frac{12(n-2)}{5}$, with equality when $n\equiv12 \pmod{20}$;\\
(\romannumeral2) (\cite{D}) for $ k\in \{3,4,5\}$, $ex_{\mathcal{P}}(n,F_k)\le \frac{13kn}{4k+2}-\frac{12k}{2k+1}$.
\end{theorem}

In Section 3, we obtain a sharp upper bound of $ex_{\mathcal{P}}(n,F_k)$,
which improves the result in Theorem \ref{103}.
Moreover, we give infinitely many extremal graphs which attain the new upper bound.

All graphs considered in this paper are undirected, finite and simple.
Let $G=(V(G),E(G))$ be a graph of vertex set $V(G)$ and edge set $E(G)$. For a vertex $v\in V(G)$, the \emph{degree} of $v$ in $G$ is denoted by $d_G(v)$. $v$ is called a \emph{$k$-vertex} of $G$ if $d_G(v)=k$. Let $\delta(G)$, $d(G)$ and $\Delta(G)$ be the minimum, average and maximum degree of $G$, respectively.
A maximal connected ($2$-connected) subgraph of $G$ is called a \emph{component} (\emph{block}) of $G$. Two graphs are \emph{vertex-disjoint} if they have no common vertex, and \emph{edge-disjoint} if they have no common edge.
For $S\subset V(G)$, the subgraph induced by $S$, denote by $G[S]$, is the graph with vertex set $S$ and edge set $\{xy\in E(G):x,y\in S\}$. Specially, write  $G\setminus v$ for the subgraph induced by $V(G)\setminus \{v\}$.
For a positive integer $t$, we use $tG$ to denote the disjoint union of $t$ copies of a graph $G$.

Let us introduce some notations on planar graphs.
 Given a plane graph $G$, the outer face of $G$ is denoted by $\Gamma(G)$. A face of size $i$ in $G$ is called an $i$-face. Let $f(G)=\sum_{i=1}f_i(G)$, where $f_i(G)$
  is the number of $i$-faces in $G$.
  Let $E_3(G)$ be the set of edges in $G$ such that each belongs to at least one $3$-face and $E_{3,3}(G)$ be the set of edges in $G$ such that each belongs to two $3$-faces.
 Furthermore, we write $E_{3,3}'(G)=E_{3,3}(G)\setminus {E(\Gamma(G))}$.
 For convenience, we denote $e(G)=|E(G)|$, $e_{3}(G)=|E_{3}(G)|$, $e_{3,3}(G)=|E_{3,3}(G)|$
   and  $e_{3,3}'(G)=|E_{3,3}'(G)|$.
 We use $n_k(G)$ to denote the number of $k$-vertices incident with exactly $k$ inner $3$-faces in a plane graph $G$.
 Given $v\in V(G)$, let $G_v$ be the graph induced by all the $3$-faces incident to $v$.

\section{ Planar Tur\'{a}n  Number of $\bm{H_k}$}

~~~~The following definition will play an important role in proving the subsequent theorem.

\begin{definition}\label{x01}
 Let $G$ be a plane graph which contains 3-faces. We recursively construct a triangular-block in the following way. Start with a $3$-face $F$ of $G$.\\
(\romannumeral1)~Take $e\in E(F)$ and search for a new $3$-face $F(e)$ of $G$ containing $e$. Add other edge(s) in this $3$-face to $E(F)$.\\
(\romannumeral2)~Repeat step (\romannumeral1), till we cannot find a new $3$-face of $G$ containing any edge in $E(F)$.
\end{definition}

The \emph{triangular-block} obtained from $F$ in Definition \ref{x01} is denoted by $\widehat{F}$.
Clearly, $\widehat{F}$ is well defined.
We shall note that an $i$-face of $\widehat{F}$ is not necessarily
an $i$-face of $G$ for $i\ge 4$, vice versa.
By Definition \ref{x01} (\romannumeral1), we have $E(\widehat{F})\subseteq E_3(G)$.
 By Definition \ref{x01} (\romannumeral2), any two triangular-blocks of $G$ are edge-disjoint.
  Let $\mathcal{B}$ be the family of  triangular-blocks of $G$. Then $e_3(G)=\sum_{\widehat{F}\in \mathcal{B}}e_3(\widehat{F})$. And if $\Gamma(G)$ is not a $3$-face,
 then we can find that $e_{3,3}(G)=\sum_{\widehat{F}\in \mathcal{B}}e_{3,3}'(\widehat{F})$.
 In fact, let $e\in E_{3,3}(G)$. Then there exists a triangular-block $\widehat{F}$ containing $e$. Since $\Gamma(G)$ is not a 3-face, there are two  inner 3-faces of $G$ containing $e$. Thus, $e\in E_{3,3}'(\widehat{F})$.

 \begin{lemma}\label{y01}
 Let $F$ be a 3-face in a plane graph $G$.
 Let $\mathcal{C}_{\widehat{F}}$ consist of $\Gamma(\widehat{F})$ and those faces of $\widehat{F}$ of size at least 4.
  If $\Gamma(G)$ is not a 3-face, then  any two faces of $\mathcal{C}_{\widehat{F}}$ are edge-disjoint.
\end{lemma}
 \noindent{\bf Proof.}
 Suppose to the contrary.
 Assume that there exist two faces of $\mathcal{C}_{\widehat{F}}$  containing a common edge $e$. Note that  $e\in E(\widehat{F})\subseteq E_3(G)$. And by the definition of ${\widehat{F}}$, there is a 3-face $U$ of $G$ containing $e$.
 Since $\Gamma(G)$ is not a 3-face, $U$ is an inner face of $G$.
By the maximality of $\widehat{F}$, $U$ is also an inner 3-face of $\widehat{F}$.
Therefore, $U\notin \mathcal{C}_{\widehat{F}}$, and the statement holds.
 ~~~~$\Box$

\begin{figure}[!ht]
	\centering
	\includegraphics[width=0.7\textwidth]{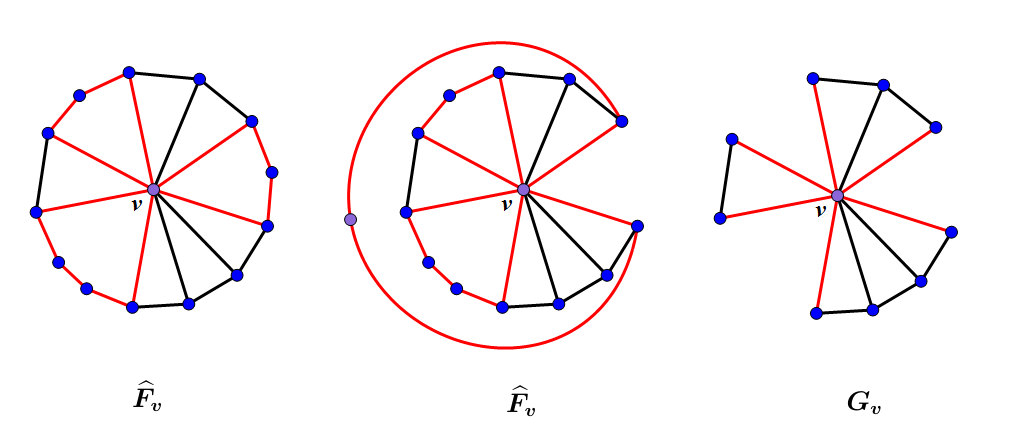}
	\caption{Two examples when $\alpha(v)=3$. }{\label{bee}}
\end{figure}

Let $\widehat{F}$ be a triangular-block of $G$ and $v\in V(\widehat{F})$ be incident with $\alpha(v)$ faces in $\mathcal{C}_{\widehat{F}}$. By Lemma \ref{y01}, if $\Gamma(G)$ is not a 3-face, then these $\alpha(v)$ faces are edge-disjoint.
 It follows that $G_v\setminus v$ contains exactly $\alpha(v)$ paths of order at least two (see Figure \ref{bee}, where $\widehat{F}_v$ is the graph induced by all the faces of $\widehat{F}$ incident to $v$).

\begin{definition}\label{x04}
 Let $\widehat{F}$ be a triangular-block in a plane graph $G$. The contribution of $\widehat{F}$ to the number of vertices in $G$, denoted by  $n(\widehat{F})$, is defined as
                $$n(\widehat{F})=\sum_{v\in V(\widehat{F})} \frac{1}{ \#{~triangular-block~sharing~v}}.$$
\end{definition}

One can observe that the number of vertices incident to 3-faces of $G$ is $\sum_{\widehat{F}\in \mathcal{B}}n(\widehat{F})$.
Thus $|V(G)|\ge \sum_{\widehat{F}\in \mathcal{B}}n(\widehat{F})$.

\begin{figure}[!ht]
	\centering
	\includegraphics[width=0.6\textwidth]{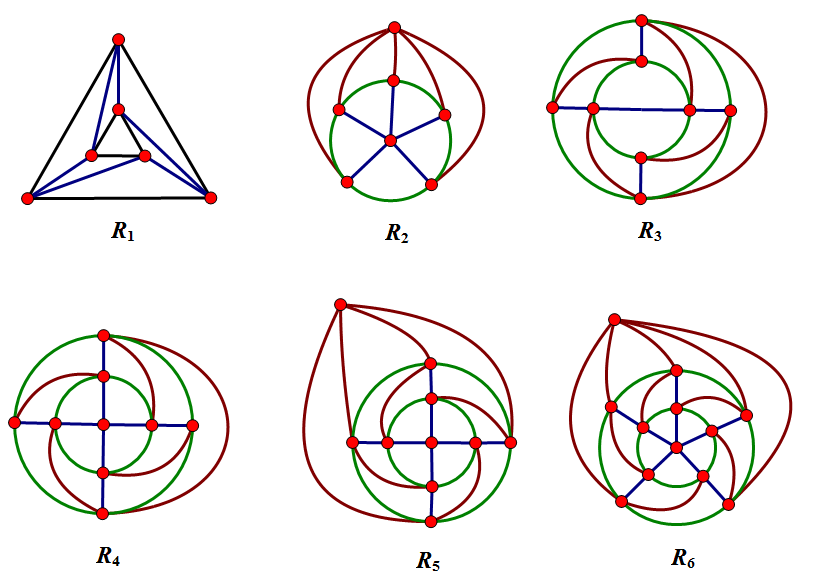}
	\caption{Extremal graphs $R_1,R_2,R_3,R_4,R_5$ and $R_6$. }{\label{ba1}}
\end{figure}

\begin{lemma}\label{107}
Let $G$ be a plane triangulation. Then $G$ is regular if and only if $G$ is isomorphic to one of $K_3,K_4,R_1$ and $R_6$ (where $R_1$ and $R_6$ see Figure \ref{ba1}).
\end{lemma}

\noindent{\bf Proof.}
Since $G$ is a plane triangulation, $e(G)=3n-6$. Assume that $G$ is $k$-regular. Thus,
 by $2e(G)=\sum_{v\in V(G)}d_G(v)=kn$, we have $n=\frac{12}{6-k}$.
Therefore, $$(n,k)\in \{(3,2),(4,3),(6,4),(12,5)\}.$$
The case when $(n,k)\in \{(3,2),(4,3)\}$ is trivial, and correspondingly $G\in \{K_3,K_4\}$.
Now we consider the case when $(n,k)=(12,5)$.
Since $2e(G)=3f_3(G)=3f(G)$, we have $f(G)=20$.
It is well known that a graph $G$ is embeddable on the plane if and only if it is embeddable on the sphere.
This implies that the outer face of $G$ can be arbitrarily chosen from its 20 faces.
Therefore,
there is a vertex $v_0\in V(G)$ such that five 3-faces incident to $v_0$ are inner faces. Moreover, $G_{v_0}\setminus v_0$ is a pentagon, say $v_1v_2v_3v_4v_5v_1$.
Similarly, there exists an inner 3-face $U_1=y_1v_1v_2y_1$ other than $v_0v_1v_2v_0$ containing $v_1v_2$.
If $y_1\in \{v_3,v_5\}$, say $y_1=v_3$, then $d_G(v_2)=3$, a contradiction. If $y_1=v_4$, then $U_1$ is an outer face, also a contradiction.
Thus, $y_1\notin N_G[v_0]$, where $N_{G}[v_0]$ is the closed neighborhood of $v_0$.
Similar discussions as above, there exists an inner 3-face $U_i=y_iv_iv_{i+1}y_i$ other than $v_0v_iv_{i+1}$ containing $v_iv_{i+1}$ for $i\in \{2,3,4,5\}$ (where $v_6=v_1$).
Moreover, $y_i\not\in N_G[v_0] \cup \{y_1,y_2,\ldots,y_{i-1}\}$ for $i\in \{2,3,4,5\}$.
Let $G'=G[N_G[v_0]\cup \{y_1,y_2,\ldots,y_{5}\}]$.
Note that $d_{G'}(v_i)=5$ for $i\in \{1,2,3,4,5\}$.
Since $G$ is a plane triangulation, we have $y_iy_{i+1}\in E(G')$ for $i\in \{1,2,3,4,5\}$.
Now $d_{G'}(y_i)=4$ for each $i$ and there is a unique vertex $z\in V(G)\setminus V(G')$.
Clearly, $zy_i\in E(G)$ for each $i$. It follows that $G\cong R_6$.

The case when $(n,k)=(6,4)$ is similar and the corresponding graph is $R_1$. Its proof is omitted here.
~~~~$\Box$

Now we give the result on the planar Tur\'{a}n number of $H_3$.
\begin{theorem}\label{106}
 $ex_{\mathcal{P}}(n,H_3)\le \frac{67n}{24}-4$ for all $n\ge 13$, with equality if and only if $n$ is divisible by $24$.
\end{theorem}

\noindent{\bf Proof.}
Let $G$ be an $H_3$-free plane graph on $n\ge 13$ vertices.
 It suffices to prove that $e(G)\le \frac{67n}{24}-4$.
 We shall proceed the proof by induction on $n$.
The statement is true when $n\in \{13,14\}$ by Theorem \ref{503}. Suppose then that $n\ge 15$.

If there exists a vertex $u\in V(G)$ with $d_G(u)\le 2$, then by the induction hypothesis,
we have $e(G\setminus u)\le \frac{67n-67}{24}-4$ and thus
$e(G)=e(G\setminus u)+d_G(u)\le \frac{67n}{24}-4$, as desired.
Recall that $e(G)\le 2n-4$ for any triangle-free graph $G$.
It follows that if $G$ contains no triangle then $e(G)<\frac{67}{24}n-4$.

If $G$ is disconnected,
let $G_1,\ldots,G_s,G_{s+1},\ldots,G_{s+t}$ be all components of $G$ such that $|V(G_i)|\le 12$ for $i\le s$ and $ |V(G_{s+j})|\ge 13$ for $j\le t$, where $s,t\ge 0$, $s+t\ge 2$ and
 $\sum_{i=1}^{s+t}|V(G_{i})|=n$.
 Since $G_i$ is a plane graph, we have $e(G_i)\le 3|V(G_i)|-6$ for all $i\le s$, and $e(G_{s+j})\le \frac{67|V(G_{s+j})|}{24}-4$ for each $j\le t$ by the induction hypothesis. Therefore,
\begin{eqnarray*}
  e(G) &\le & 3\sum_{i=1}^{s}|V(G_{i})|-6s+
\frac{67\sum_{j=1}^{t}|V(G_{s+j})|}{24}-4t \\
  &=& (\frac{67n}{24}-4t-2s)+(\frac{5\sum_{i=1}^{s}|V(G_{i})|}{24}-4s) \\
  &\le& (\frac{67n}{24}-4)+(\frac{5}{24}\times 12s-4s)\\
  &\le& \frac{67n}{24}-4.
\end{eqnarray*}

In the following we may assume that  $G$ is a connected graph with $\delta(G)\ge 3$ and $f_3(G)>0$.
If $G$ is a plane triangulation, then $e(G)=3n-6$ and hence
there is a vertex $v$ such that $d_G(v)\ge 6$.
This implies that $G$ contains a copy of $H_3$ with center $v$, a contradiction.
Therefore $G$ is not a plane triangulation.
Without loss of generality, we may  assume that $\Gamma(G)$ is not a $3$-face
and $F$ is an inner 3-face of $G$.
Assume that the size of $\Gamma(\widehat{F})$ is $l_0$ and the sizes of other faces in $\mathcal{C}_{\widehat{F}}$ are  $l_1,l_2,\ldots,l_{|\mathcal{C}_{\widehat{F}}|-1}$, respectively.
 Let $v$ be a vertex incident with exactly $\alpha(v)$ faces in $\mathcal{C}_{\widehat{F}}$.
 Recall that $G_v\setminus v$ contains exactly $\alpha(v)$ paths of order at least two.
 Thus, $\alpha(v)\le 2$ for any $v\in V(\widehat{F})$ since $G$ is $H_3$-free.

\noindent{\bf {Claim 1.}} Let $A_{\widehat{F}}=\{v\in V(\widehat{F}):\alpha(v)=2\}$.
Then $4\le d_{\widehat{F}}(v)\le 6$ for any $v\in A_{\widehat{F}}$ and $d_{\widehat{F}}(v)\le 5$ for any $v\in V(\widehat{F})\setminus A_{\widehat{F}}$.

\noindent{\bf Proof.}
Clearly, $d_{\widehat{F}}(v)\ge 4$ for any $v\in A_{\widehat{F}}$.
If $d_{\widehat{F}}(v)\ge 7$ for some $v\in A_{\widehat{F}}$,
then $G_v\setminus v$ contains a copy of $P_5\cup P_2$ or $P_4\cup P_3$. It follows that $G_v$ contains a copy of $H_3$, a contradiction.
Thus, $4\le d_{\widehat{F}}(v)\le 6$ for any $v\in A_{\widehat{F}}$.
On the other hand,
if $d_{\widehat{F}}(v)\ge 6$ for some $v\in V(\widehat{F})\setminus A_{\widehat{F}}$,
then $G_v\setminus v$ contains a copy of $P_6$. It follows that $G_v$ contains a copy of $H_3$, a contradiction.
Thus, $d_{\widehat{F}}(v)\le 5$ for any $v\in V(\widehat{F})\setminus A_{\widehat{F}}$.
~~~~$\Box$

\begin{figure}[!ht]
	\centering
	\includegraphics[width=0.5\textwidth]{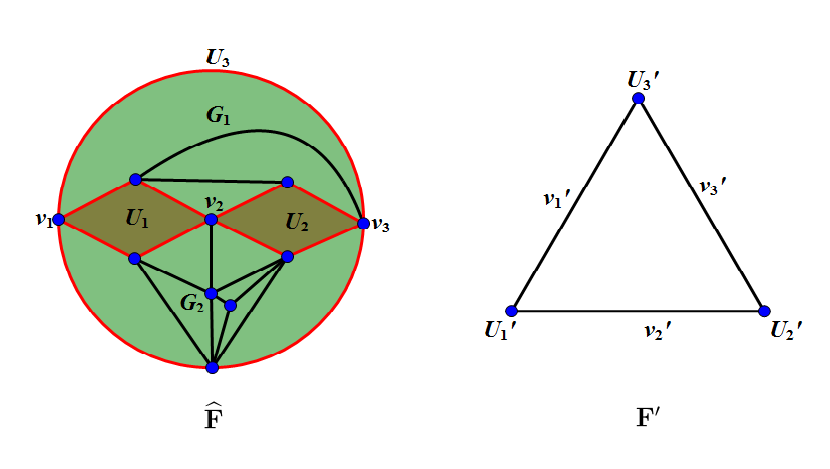}
	\caption{ An example of graphs $\widehat{F}$ and $F'$. }{\label{bah}}
\end{figure}

\noindent{\bf {Claim 2.}} $|A_{\widehat{F}}|\le |\mathcal{C}_{\widehat{F}}|-1$.

\noindent{\bf Proof.}
We define a graph $F'$ as follows (see Figure \ref{bah}): corresponding to each face $U_i$ in $\mathcal{C}_{\widehat{F}}$
there is a vertex $U_{i}'$ of $F'$, and corresponding to each vertex $v_k$ of $A_{\widehat{F}}$ there is an edge $v_{k}'$ of $F'$; two vertices $U_{i}'$ and $U_{j}'$ are joined by an edge $v_{k}'$ in $F'$ if and only if their corresponding faces $U_i$ and $U_j$ share a vertex $v_k\in V(\widehat{F})$.
In order to prove the claim, it suffices to  show $F'$ is acyclic.

Suppose to the contrary that $F'$ contains a cycle. Then $\mathcal{C}_{\widehat{F}}$ separates ${\widehat{F}}$ into at least two edge-disjoint subgraphs $G_1$ and $G_2$ (see Figure \ref{bah}).
Let $e$ be an arbitrary edge of some face in $\mathcal{C}_{\widehat{F}}$.
Thus, $e\in E_{3}(G)$.
Without loss of generality, we may assume that $e\in E(G_1)$.
Note that $\Gamma(G)$ is not a $3$-face.
By the discussion in Lemma \ref{y01}, each inner 3-face of $G$ containing $e$ is also an inner 3-face of $\widehat{F}$. According to the choice of $e$, there is exactly one 3-face $F(e)$ of $G$ containing $e$, that is, the other face of $G$ containing $e$ is not a 3-face.
This indicates that $G_1$ is a triangular-block, which contradicts the maximality of $\widehat{F}$.
~~~~$\Box$

\noindent{\bf {Claim 3.}} $n_5(G)\le \frac{3}{4}n$, with equality if and only if
 any triangular-block of $G$ is isomorphic to $R_6$ (see Figure \ref{ba1}) and
any vertex of $G$ belongs to exactly one triangular-block of $G$.

\noindent{\bf Proof.}
Recall that  any two triangular-blocks of $G$ are edge-disjoint.
For any vertex $v\in V(G)$, $v$ is incident with at most two triangular-blocks since $G$ is $H_3$-free.
Let $B_{\widehat{F}}=\{v\in V({\widehat{F}}):$
$v$ is shared by exactly two triangular-blocks of $G\}$.
 Let $v\in B_{\widehat{F}}$. Then $v$ lies in $\Gamma(\widehat{F})$. Let $e_1,e_2$ be the two edges of $\Gamma(\widehat{F})$ incident to $v$.
 Then there is an inner 3-face $U_i$ containing $e_i$ for $i\in \{1,2\}$.
  If $d_{\widehat{F}}(v)\ge 4$, then $U_1$ and $U_2$ are edge-disjoint.
Note that the other triangular-block containing $v$ also has an inner 3-face $U_3$ containing $v$.
 Therefore, $G$ contains a copy of $H_3$, a contradiction.
 It follows that $d_{\widehat{F}}(v)\le 3$ for any vertex $v\in B_{\widehat{F}}$.
 Since $G$ is $H_3$-free.  Moreover, $A_{\widehat{F}}\cap B_{\widehat{F}}=\phi$, since by Claim 1 $d_{\widehat{F}}(v)\ge 4$ for any vertex $v\in A_{\widehat{F}}$.
Note that for every $i\in \{0,1,2,\ldots,|\mathcal{C}_{\widehat{F}}|-1\}$,  we can add $l_i-3$ edges to the face of size $l_i$ such that the resulting graph is a plane triangulation. Therefore,
$$e(\widehat{F}) =(3|V(\widehat{F})|-6)-\sum_{i=0}^{|\mathcal{C}_{\widehat{F}}|-1}(l_i-3).  $$
On the other hand, by Claim 1 we have,
$$2e(\widehat{F})=\sum_{v\in V(\widehat{F})}d_{\widehat{F}}(v)\le 6|A_{\widehat{F}}|+3|B_{\widehat{F}}|+5(|V(\widehat{F})|-|A_{\widehat{F}}|-|B_{\widehat{F}}|)
=5|V(\widehat{F})|+|A_{\widehat{F}}|-2|B_{\widehat{F}}|.$$
Combining with above two inequalities, we have
\begin{align}\label{dda}
  |V(\widehat{F})|\le 2l_0+2\sum_{i=1}^{|\mathcal{C}_{\widehat{F}}|-1}l_i
  -6|\mathcal{C}_{\widehat{F}}|+|A_{\widehat{F}}|-2|B_{\widehat{F}}|+12.
\end{align}
According to the definition of $\widehat{F}$, we can find that all faces of $\widehat{F}$ are cycles. Therefore,
by (\ref{dda}) and Claim 2, we have
\begin{eqnarray}\label{ddb}
  n_5(\widehat{F}) &\le & |V(\widehat{F})|-\Big(
\sum_{i=0}^{|\mathcal{C}_{\widehat{F}}|-1}l_i-|A_{\widehat{F}}|\Big) \nonumber\\
  &=&  \frac{3}{4}\Big({|V(\widehat{F})|-\frac{|B_{\widehat{F}}|}{2}}\Big)
  +\frac{1}{4}\Big\{ |V(\widehat{F})|-4l_0-4\sum_{i=1}^{|\mathcal{C}_{\widehat{F}}|-1}{l_i}
    +4|A_{\widehat{F}}|+\frac{3|B_{\widehat{F}}|}{2}\Big\}      \nonumber  \\
  &\le&  \frac{3}{4}\Big({|V(\widehat{F})|-\frac{|B_{\widehat{F}}|}{2}}\Big)-\frac{1}{4}\Big\{ (2l_0-6)+(6|\mathcal{C}_{\widehat{F}}|-6-5|A_{\widehat{F}}|)+\frac{1}{2}|B_{\widehat{F}}|
+2\sum_{i=1}^{|\mathcal{C}_{\widehat{F}}|-1}l_i\Big\} \nonumber \\
  &\le& \frac{3}{4}\Big({|V(\widehat{F})|-\frac{|B_{\widehat{F}}|}{2}}\Big).
\end{eqnarray}
The last inequality holds in equality if and only if $l_0=3$ and $|A_{\widehat{F}}|=|B_{\widehat{F}}|=|\mathcal{C}_{\widehat{F}}|-1=0$.
 If  ${n_5(\widehat{F})}= \frac{3}{4}({|V(\widehat{F})|-\frac{|B_{\widehat{F}}|}{2}})$,
 then (1) becomes an equality, and thus $|V(\widehat{F})|=12$.
 Furthermore, $l_0=3$ and $|\mathcal{C}_{\widehat{F}}|=1$ imply that $\widehat{F}$ is a plane triangulation,
 that is, $e(\widehat{F})=3|V(\widehat{F})|-6=30$.
 Recall that the degree of every vertex in  $V(\widehat{F})\setminus A_{\widehat{F}}$ is at most 5.
 Thus, $\widehat{F}$ is 5-regular. By Lemma \ref{107}, we have $\widehat{F}\cong R_6$ for any $\widehat{F}\in \mathcal{B}$.

Let $v$ be a 5-vertex incident with five inner 3-face of $G$. Clearly, $v$ belongs to only one triangular-block of $G$.  Thus
$n_5(G)=\sum_{\widehat{F}\in \mathcal{B}}n_5(\widehat{F})$.
 By Definition \ref{x04} and the choice of $B_{\widehat{F}}$, we have  $n(\widehat{F})=|V(\widehat{F})|-\frac{|B_{\widehat{F}}|}{2}$.
It follows that
 \begin{eqnarray}\label{ddc}
 n\ge \sum_{\widehat{F}\in \mathcal{B}}n(\widehat{F})=\sum_{\widehat{F}\in \mathcal{B}}(|V(\widehat{F})|-\frac{|B_{\widehat{F}}|}{2}).
 \end{eqnarray}
 Thus by (\ref{ddb}), we have
 \begin{align}\label{ddd}
   n_5(G)=\sum_{\widehat{F}\in \mathcal{B}}n_5(\widehat{F})\le \frac{3}{4}\sum_{\widehat{F}\in \mathcal{B}}(|V(\widehat{F})|-\frac{|B_{\widehat{F}}|}{2})
 \le \frac{3}{4}n.
 \end{align}
 ~~~~We now characterize the equality in (\ref{ddd}).
 ${n_5(G)}= \frac{3}{4}n$ if and only if equalities hold in (\ref{ddc}) and (\ref{ddb}) for every $\widehat{F}\in \mathcal{B}$.
Note that if (\ref{ddb}) becomes an equality,
then $\widehat{F}\cong R_6$ and $B_{\widehat{F}}=0$.
 Moreover, (\ref{ddc}) becomes an equality means that $n=\sum_{\widehat{F}\in \mathcal{B}}|V(\widehat{F})|$, i.e.,
  any vertex of $G$ belongs to exactly one triangular-block of $G$.
  Conversely, assume that $\widehat{F}\cong R_6$ for each $\widehat{F}\in \mathcal{B}$ and
  $v$ belongs to exactly one $\widehat{F}$ for each $v\in V(G)$.
Then
   $$n_5(G)=\sum_{\widehat{F}\in \widehat{B}}n_5(\widehat{F})=\sum_{\widehat{F}\in \widehat{B}}\frac{3}{4}|V(\widehat{F})|=\frac{3}{4}n,$$
   since $|V(R_6)|=12$ and $n_5(R_6)=9$.
~~~~$\Box$

Note that for $i\le 5$, each $i$-vertex is incident with at most $i$ inner 3-faces of $G$. Moreover, for all $i\ge 6$, each $i$-vertex is incident with at most four 3-faces (otherwise, if there is a vertex $v$ such that $d_G(v)\ge 6$ and $v$ is incident with at least five 3-faces, then $G_v$ contains a copy of $H_3$). Double counting the number of vertices, we have
$$3f_3(G)\le 5n_5(G)+4(n-n_5(G))= 4n+n_5(G). \nonumber$$
By Claim 3, we have $3f_3(G)\le \frac{19n}{4}$, that is, $f_3(G)\le\frac{19n}{12}$.
Moreover,
\begin{align}\label{dde}
  2e(G)=3f_3(G)+\sum_{i\ge 4}if_i(G)\ge 3f_3(G)+4(f(G)-f_3(G))=4f(G)-f_3(G).
\end{align}
This implies that $f(G)\le \frac{2e(G)+f_3(G)}{4}$.

By Euler's formula, $n-2=e(G)-f(G)\ge e(G)-\frac{2e(G)+\frac{19n}{12}}{4}$.
This indicates that $e(G)\le \frac{67n}{24}-4$, as desired.

\begin{figure}[!ht]
	\centering
	\includegraphics[width=0.5\textwidth]{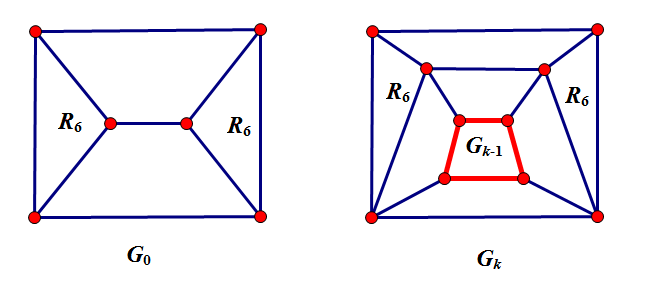}
	\caption{The construction of $G_k$. }{\label{bdc}}
\end{figure}

From the proof above, we see that equality in $e(G)\le \frac{67n}{24}-4$ is achieved if and only if $n_5(G)= \frac{3}{4}n$ and equality holds in (\ref{dde}). This implies that $e(G)= \frac{67n}{24}-4$
if and only if $G$ is a connected $H_3$-free plane graph satisfying:  any triangular-block of $G$ is isomorphic to $R_6$;
any vertex of $G$ belongs to exactly one triangular-block of $G$;
any face of $G$ is either a 3-face or a 4-face.
We next construct such an extremal plane graph $G$.
 Let $G_0$ be the plane graph obtained from two disjoint copies of $R_6$ by adding three independent edges between their outer faces, as depicted in Figure \ref{bdc}. We construct $G_k$ of order $n$ recursively for all $k\ge1$ via the illustration given in Figure \ref{bdc}: the entire graph $G_{k-1}$ is placed into the center bold red quadrangle of Figure \ref{bdc} (in such a way that the center bold red quadrangle is identified with the outer quadrangle of $G_{k-1}$). One can check that $G_k$ is $H_3$-free with $n=24(k+1)$ vertices and $\frac{67n}{24}-4$ edges for all $k\ge 0$.
 This completes the proof.~~~~$\Box$

\section{ Planar Tur\'{a}n  Number of $\bm{F_k}$}
~~~~In this section we will determine the planar Tur\'{a}n number of $F_k$.
We first need to introduce a lemma.

\begin{lemma}\label{101} (\cite{bab})
There does not exist a planar graph on $n\in \{11,13\}$ vertices with exactly one vertex of degree $4$ and $n-1$ vertices of degree $5$.
\end{lemma}

Recall that $ex_{\mathcal{P}}(n,F_k)=3n-6$ for all $k\ge 6$.
The following theorem characterizes the case when $ex_{\mathcal{P}}(n,F_k)=3n-6$
for $k\le 5$.

\begin{theorem}\label{303}
Let $n,k$ be two positive integers with $k\leq5$ and $n\ge k+2$. Then
$ex_{\mathcal{P}}(n,F_k)=3n-6$ if and only if $(n,k)\in \{(6,4),(7,5),(8,5),(9,5),(10,5),(12,5)\}$.
\end{theorem}

\noindent{\bf Proof.}
Note that if $G$ is an $F_k$-free plane triangulation on $n$ vertices, then $G$ contains a copy of $F_{\Delta(G)-1}$ and hence $\Delta(G)\le k$.
This, together with $\Delta(G)\ge d(G)=6-\frac{12}{n}$, implies that $n\le \frac{12}{6-k}$.
Thus, $(n,k)\in \{(6,4),(7,5),(8,5),(9,5),(10,5),(11,5),(12,5)\}$.
The plane triangulation $R_1$ depicted in Figure \ref{ba1} is $F_4$-free. Thus
$ex_{\mathcal{P}}(6,F_4)=12$.
For $n\in \{7,8,9,10,12\}$, the plane triangulations $R_2$-$R_6$ depicted in Figure \ref{ba1} are $F_5$-free and thus $ex_{\mathcal{P}}(n,F_5)=3n-6$.
We next show that  $ex_{\mathcal{P}}(11,F_5)<3n-6$.
If $ex_{\mathcal{P}}(11,F_5)=3n-6=27$, then $d(G)=6-\frac{12}{n}$.
Thus, $\Delta(G)\ge 5$ and $\delta(G)\le 4$.
Since $\Delta(G)\le k=5$, we have $\Delta(G)=5$.
Therefore,
$$54=2e(G)=\sum_{v\in V(G)}d_G(v)\le \delta(G)+(n-1)\Delta(G)\le 4+5(n-1)=54.$$
 This implies that $G$ is a plane triangulation with exactly one vertex of degree $4$ and $n-1$ vertices of degree $5$,
contradicting to Lemma \ref{101}. ~~~~$\Box$

In the following, we introduce a new definition of blocks in a plane graph.

\begin{definition}\label{x03}
Let $G$ be a plane graph  containing 3-faces and $\Gamma(G)$ is not a 3-face. We recursively construct an improvement-block in the following way. \\
(\romannumeral1)~Start with a $3$-face $F$ of $G$.\\
(\romannumeral2)~Take $v\in V(\widehat{F})$ such that $v$ is incident with exactly $l\ge 2$ faces in $\mathcal{C}_{\widehat{F}}$. Recall that $G_v\setminus v$ consists of $l$ path-components $P^1,P^2,\ldots,P^l$.
Then delete $v$, add $l$ isolated vertices $v_1,v_2,\ldots,v_l$ and link $v_i$ to each vertex of $P^i$ for each $i\in \{1,2,\ldots,l\}$.\\
(\romannumeral3)~Repeat step (\romannumeral1), till the resulting graph contains no vertex satisfying the above property.
\end{definition}

\begin{figure}[!ht]
	\centering
	\includegraphics[width=0.5\textwidth]{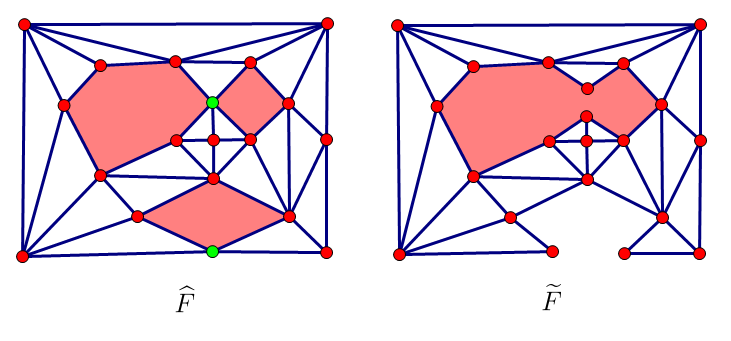}
	\caption{A triangular-block $\widehat{F}$ and corresponding improvement-block $\widetilde{F}$. }{\label{bae}}
\end{figure}

The \emph{improvement-block} obtained from $F$ in Definition \ref{x03} is denoted by $\widetilde{F}$.
Clearly, $\widetilde{F}$ is well defined.
For instance,  a triangular-block and corresponding improvement-block of a plane graph $G$ are shown in Figure \ref{bae}. We can observe that $e_{3,3}'(\widehat{F})=e_{3,3}'(\widetilde{F})$ and $e_{3}(\widehat{F})=e_{3}(\widetilde{F})$.
 Moreover, any two faces of $\mathcal{C}_{\widetilde{F}}$ are vertex-disjoint by Definition \ref{x03} (\romannumeral2),
where $\mathcal{C}_{\widetilde{F}}$ consists of $\Gamma(\widetilde{F})$ and those faces of $\widetilde{F}$ of size at least 4.

Now we give our main theorem of this section.
We note that for $k=2$, the bound in Theorem \ref{103} meets the bound in the following theorem.

\begin{theorem}\label{102}
Let $n,k$ be integers with $k\in \{2,3,4,5\}$ and $n\ge \frac{12}{6-k}+1$. Then
$ ex_{\mathcal{P}}(n,F_k)\le \frac{24k}{7k+6}(n-2)$,
with equality if $n\equiv{\frac{12(k+2)}{6-k}}\pmod{{\frac{28k+24}{6-k}}}$.
\end{theorem}

\noindent{\bf Proof.}
Let $G$ be an $F_k$-free plane graph on $n\ge \frac{12}{6-k}+1$ vertices, where $k\in \{2,3,4,5\}$.
By the proof of Theorem \ref{303}, $G$ is not a plane triangulation.
We shall proceed the proof by induction on $n$.
Assume first $n=\frac{12}{6-k}+1$. If $k=2$, then it is easy to verify that $n=4$ and $e(G)\le 4$, as desired. If $k\in \{3,4,5\}$, then
 $$e(G)\le 3(\frac{12}{6-k}+1)-7=\frac{4k+12}{6-k}\le \frac{24k}{7k+6}(n-2).$$

Suppose then that $n\ge \frac{12}{6-k}+2$.
If there exists a vertex $u\in V(G)$ with $d_G(u)\le 2$, then by the induction hypothesis, $e(G\setminus u)\le \frac{24k}{7k+6}(n-3)$ and thus $e(G)=e(G\setminus u)+d_G(u)\le \frac{24k}{7k+6}(n-2)$, as desired.
Next assume that $\delta(G)\ge 3$.
If $G$ is disconnected,
let $G_1,\ldots,G_s,G_{s+1},\ldots,G_{s+t}$ be all components of $G$ such that
$|V(G_i)|\le \frac{12}{6-k}$ for $i\le s$ and $ |V(G_{s+j})|\ge \frac{12}{6-k}+1$ for $j\le t$, where $s,t\ge 0$, $s+t\ge 2$ and $\sum_{i=1}^{s+t}|V(G_{i})|=n$.
 Since $G_i$ is a plane graph, we have $e(G_i)\le 3|V(G_i)|-6$ for all $i\le s$, and $e(G_{s+j})\le \frac{24k}{7k+6}(|V(G_j)|-2)$ for all $j\le t$ by the induction hypothesis. Therefore,
\begin{eqnarray*}
 e(G) &=& 3\sum_{i=1}^{s}|V(G_{i})|-6s+
\frac{24k(\sum_{j=1}^{t}|V(G_{s+j})|-2t)}{7k+6} \\
    &=&
    \frac{24k(n-2)}{7k+6}+(\frac{18-3k}{7k+6}\sum_{i=1}^{s}|V(G_{i})|-6s)-\frac{48k}{7k+6}(t-1)\\
    &\le&   \frac{24k(n-2)}{7k+6}-\frac{42k}{7k+6}s-\frac{48k}{7k+6}(t-1)\\
    &\le&   \frac{24k(n-2)}{7k+6}-\frac{42k}{7k+6}(s+t)+\frac{48k}{7k+6}\\
    &<& \frac{24k(n-2)}{7k+6}.
\end{eqnarray*}
~~~~It remains to consider the case when $G$ is connected and $\delta(G)\ge 3$.
Note that $G$ is not a plane triangulation.
Without loss of generality, we may  assume that $\Gamma(G)$ is not a $3$-face.
If $G$ contains no 3-face, then $e(G)\le 2n-4<\frac{24k}{7k+6}(n-2)$, as desired.
Now assume that $G$ contains a 3-face $F$.
Note that $G$ is $F_k$-free. Then both $\widehat{F}$ and $\widetilde{F}$ are $F_k$-free.

\noindent{\bf {Claim 1.}} $e_{3,3}'(\widehat{F})\le \frac{3k-6}{2k}e_{3}(\widehat{F})$ with equality if and only if $\widehat{F}\cong J_k$, where $J_k$ is isomorphic to
$K_3$ for $k=2$,
$K_4$ for $k=3$, $R_1$ for $k=4$ and $R_6$ for $k=5$.

\noindent{\bf Proof.}
Equivalently, we shall prove that $e_{3,3}(\widetilde{F})\le \frac{3k-6}{2k}e_{3}(\widetilde{F})$.
 Recall that any two faces of $\mathcal{C}_{\widetilde{F}}$ are vertex-disjoint.
 Thus, for any $j$-vertex $v$ of $\widetilde{F}$,
 there is at most one face of $\mathcal{C}_{\widetilde{F}}$ incident to $v$.
This indicates that other $j-1$ faces of $\widetilde{F}$ incident to $v$ 
are consecutive 3-faces.
It follows that $\widetilde{F}$ contains a copy of the $(j-1)$-fan $F_{j-1}$.
 Since $\widetilde{F}$ is $F_k$-free, we have $j-1<k$.
 It follows that $\Delta(\widetilde{F})\le k$ and $2e_3(\widetilde{F})=2e(\widetilde{F})\le k|V(\widetilde{F})|$.
Therefore, $e_3(\widetilde{F})\le \frac{k}{2}|V(\widetilde{F})|$.
Assume that the size of $\Gamma(\widetilde{F})$ is $l_0$ and the sizes of other faces of $\widetilde{F}$ are  $l_1,l_2,\ldots,l_{|\mathcal{C}_{\widetilde{F}}|-1}$, respectively.
Note that $e_3(\widetilde{F})=(3|V(\widetilde{F})|-6)-\sum_{i=0}^{|\mathcal{C}_{\widetilde{F}}|-1}(l_i-3)$.
Combining with $e_3(\widetilde{F})\le \frac{k}{2}|V(\widetilde{F})|$,
we have 
\begin{align}\label{ddn}
  \sum_{i=0}^{|\mathcal{C}_{\widetilde{F}}|-1}l_i\ge \frac{6-k}{2}|V(\widetilde{F})|+3|\mathcal{C}_{\widetilde{F}}|-6.
\end{align}
Recall that all faces of $\widehat{F}$ are cycles.
Clearly, $e_{3}(\widehat{F})-e_{3,3}'(\widehat{F})=\sum_{j=0}^{|\mathcal{C}_{\widehat{F}}|-1}l_j'
=\sum_{i=0}^{|\mathcal{C}_{\widetilde{F}}|-1}l_i$,
where the size of $\Gamma(\widehat{F})$ is $l_0'$ and the sizes of other faces of $\widehat{F}$ are  $l_1',l_2',\ldots,l_{|\mathcal{C}_{\widehat{F}}|-1}'$, respectively.
It follows from (\ref{ddn}) that
\begin{align}\label{ddx}
  e_{3}(\widetilde{F})-e_{3,3}'(\widetilde{F})=\sum_{i=0}^{s}l_i\ge \frac{6-k}{2}|V(\widetilde{F})|+3|\mathcal{C}_{\widetilde{F}}|-6.
\end{align}
If $|V(\widetilde{F})|> \frac{12}{6-k}$, then
$$\frac{e_{3}(\widetilde{F})-e_{3,3}'(\widetilde{F})}{e_{3}(\widetilde{F})}
    \ge \frac{\frac{6-k}{2}|V(\widetilde{F})|-3}{\frac{k}{2}|V(\widetilde{F})|}
    >\frac{6-k}{k}-\frac{3}{\frac{k}{2}\frac{12}{6-k}}
    =\frac{6-k}{2k},$$
that is, $e_{3,3}'(\widetilde{F})< \frac{3k-6}{2k}e_{3}(\widetilde{F})$, as desired.
Note that $\sum_{i=0}^{|\mathcal{C}_{\widetilde{F}}|-1}l_i\ge l_0\ge3$. If $|V(\widetilde{F})|\le \frac{12}{6-k}$, then
$$\frac{e_{3}(\widetilde{F})-e_{3,3}'(\widetilde{F})}{e_{3}(\widetilde{F})}
=\frac{\sum_{i=0}^{|\mathcal{C}_{\widetilde{F}}|-1}l_i}{e_{3}(\widetilde{F})}
    \ge \frac{3}{ \frac{k}{2}|V(\widetilde{F})|}\ge \frac{6-k}{2k}.$$
 Therefore,
$e_{3,3}'(\widetilde{F})\le \frac{3k-6}{2k}e_{3}(\widetilde{F})$, with equality if and only if $\sum_{i=0}^{|\mathcal{C}_{\widetilde{F}}|-1}l_i=l_0=3$, $|V(\widetilde{F})|=\frac{12}{6-k}$ and $e(\widetilde{F})=e_{3}(\widetilde{F})=\frac{k}{2}|V(\widetilde{F})|$.
This indicates that $\widetilde{F}$ is a $k$-regular plane triangulation.
 By Lemma \ref{107}, we have $\widetilde{F}\cong J_k$. Furthermore, we can see that $\widehat{F}\cong J_k$. Conversely, assume that $\widehat{F}\cong J_k$. It is easy to check that $e_{3,3}'(\widehat{F})= \frac{3k-6}{2k}e_{3}(\widehat{F})$.
~~~~$\Box$

Recall that $e_3(G)=\sum_{\widehat{F}\in \mathcal{B}}e_3(\widehat{F})$ and $e_{3,3}(G)=\sum_{\widehat{F}\in \mathcal{B}}e_{3,3}'(\widehat{F})$ since $\Gamma(G)$ is not a 3-face.
By Claim 1, we have
 \begin{align}\label{ddh}
  e_{3,3}(G)=\sum_{\widehat{F}\in \mathcal{B}}e_{3,3}'(\widehat{F})
  \le \frac{3k-6}{2k}\sum_{\widehat{F}\in \mathcal{B}}e_{3}(\widehat{F})
  =\frac{3k-6}{2k}e_3(G). \nonumber
 \end{align}
Hence, $3f_3(G)=e_3(G)+e_{3,3}(G)\le \frac{5k-6}{2k}e_3(G)$, that is, $f_3(G)\le \frac{5k-6}{6k}e_3(G)$.
Note that $e_3(G)\le e(G)$.
By (\ref{dde}) we have
$$f(G)\le \frac{2e(G)+f_3(G)}{4}\le \frac{2e(G)+\frac{5k-6}{6k}e_3(G)}{4}\le \frac{(17k-6)e(G)}{24k}.$$
By Euler's formula, $n-2=e(G)-f(G)\ge \frac{(7k+6)e(G)}{24k}$. Thus,
$e(G)\le \frac{24k(n-2)}{7k+6}$, as desired.

\begin{figure}[!ht]
	\centering
	\includegraphics[width=0.5\textwidth]{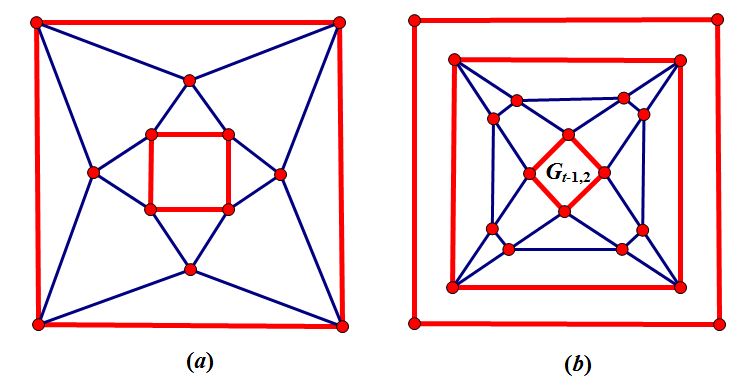}
	\caption{Construction of $G_{t,2}$. }{\label{ba2}}
\end{figure}

From the proof above, we see that $e(G)= \frac{24k(n-2)}{7k+6}$ if and only if  $e_3(G)= e(G)$, $f(G)= \frac{2e(G)+f_3(G)}{4}$ and $e_{3,3}'(\widehat{F})= \frac{3k-6}{2k}e_{3}(\widehat{F})$ for any $\widehat{F}\in \mathcal{B}$. This implies that if $e(G)= \frac{24k(n-2)}{7k+6}$,
then $G$ is a connected $F_k$-free plane graph satisfying:
$\widehat{F}\cong J_k$ for any $\widehat{F}\in \mathcal{B}$; each edge in $G$ belongs to either $E_{3,3}(G)$ or $E_{3,4}(G)$.
The construction of extremal plane graphs for $k=2$ can be found  in ~\cite{BT}.
We now follow this direction to construct the extremal graphs for general $k$.

 Let $G_{0,2}$ be the graph depicted in Figure \ref{ba2} (a).
 Let us then proceed to define further plane graphs $G_{t,2}$ on $20t+12$ vertices inductively via
 the illustration given in Figure \ref{ba2} (b). Here,
 the entire graph $G_{t-1,2}$ is placed into the center bold red quadrangle of Figure \ref {ba2} (b), and the entire graph $G_{0,2}$ is then placed between the two given bold red quadrangles of
  Figure \ref{ba2} (b) (in such a way that these are identified with the bold red quadrangles).
One can observe that $G_{t,2}$ is an $F_2$-free plane graph of order $20t+12$ and size $48t+24$.
We then construct $G_{t,k}$ by replacing every 3-face of $G_{t,2}$ with a copy of $J_k$. One can observe that $G_{t,k}$ is an $F_k$-free plane graph of order $\frac{28k+24}{6-k}t+\frac{12(k+2)}{6-k}$ and size $\frac{96k}{6-k}t+\frac{48k}{6-k}$ for $k\in \{2,3,4,5\}$.
Clearly, $e(G_{t,k})=\frac{24k}{7k+6}(|V(G_{t,k})|-2)$.

This completes the proof.
 ~~~~$\Box$

\small {

}
\end{document}